\documentclass[11pt,twoside, a4paper, english, reqno]{amsart}
\usepackage{graphicx, amsmath, varioref, amscd, amssymb,color, bm, stmaryrd, amsthm, epsfig, color, epic, enumerate}
\usepackage{fancybox}
\usepackage[pdftex, colorlinks=true]{hyperref}
\usepackage{upref}
\usepackage{pdfsync}

\usepackage{eso-pic}
\usepackage{color}
\usepackage{type1cm}
\makeatletter
  \AddToShipoutPicture*{%
    \setlength{\@tempdimb}{.12\paperwidth}%
    \setlength{\@tempdimc}{.5\paperheight}%
    \setlength{\unitlength}{1pt}%
    \put(\strip@pt\@tempdimb,\strip@pt\@tempdimc){%
      \makebox(0,0){\rotatebox{90}{\textcolor[gray]{0.75}{\fontsize{2cm}{2cm}}}}
    }
}

\makeatother

\numberwithin{equation}{section}
\allowdisplaybreaks

\newtheorem{theorem}{Theorem}[section]
\newtheorem{lemma}[theorem]{Lemma}

\newtheorem{proposition}[theorem]{Proposition}
\theoremstyle{definition}
\newtheorem{definition}[theorem]{Definition}

\theoremstyle{remark}

\newcommand{\Div}{\operatorname{div}}

\newcommand{\vr}{\varrho}
\newcommand{\vu}{\vc{u}}

\newcommand{\vc}[1]{{\bm{#1}}}

\newcommand{\R}{\mathbb{R}}
\newcommand{\N}{\mathbb{N}}

\begin{document}

\title[On a Free Boundary Problem for a Model of Polymeric Fluids]{On  a free boundary problem for polymeric fluids: Global existence of weak solutions }

\author[Donatelli]{Donatella Donatelli}
\address[Donatelli]{\newline
Departement of Engineering Computer Science and Mathematics\\
University of L'Aquila\\
67100 L'Aquila, Italy.}
\email[]{\href{donatell@univaq.it}{donatell@univaq.it}}
\urladdr{\href{http://univaq.it/~donatell}{univaq.it/\~{}donatell}}

\author[Trivisa]{Konstantina Trivisa}
\address[Trivisa]{\newline
Department of Mathematics \\ University of Maryland \\ College Park, MD 20742-4015, USA.}
\email[]{\href{http://www.math.umd.edu}{trivisa@math.umd.edu}}
\urladdr{\href{http://www.math.umd.edu/~trivisa}{math.umd.edu/\~{}trivisa}}

\date{\today}

\subjclass[2010]{Primary: 35Q30, 76N10; Secondary: 46E35.}

\keywords{Doi model; suspensions of rod-like molecules; fluid-particle interaction model; compressible Navier-Stokes equations; Fokker-Planck-type equation, free boundary problems.}

\thanks{}

\maketitle

\begin{abstract}
We investigate the stability and global existence of weak solutions to a free boundary problem governing the evolution of  polymeric fluids. We construct weak solutions of  the two-phase  model by performing the asymptotic limit of  a macroscopic model governing  the suspensions of rod-like molecules  (known as Doi-Model) in compressible fluids as the adiabatic exponent $\gamma$ goes to $\infty.$  The convergence of these solutions, up to a subsequence, to the free-boundary problem is established using techniques in the spirit of Lions and Masmoudi
\cite{LionsMasmoudi-1999}.

\end{abstract}

\tableofcontents{}

\section{Introduction}\label{S1}

\subsection{Motivation}
The evolution of rod-like molecules in polymeric fluids is of great scientific interest with a variety of applications in science and engineering.  The present article deals with a free-boundary problem for the suspension of rod-like molecules in a dilute regime. 
This article is part of a research program whose objective is the investigation of fluids with imbedded domains (large bubbles) filled with gas in the presence of rod-like molecules: standard models involve a threshold on the pressure beyond which one has the incompressible Navier-Stokes equations for the fluid and below which one has a compressible model for the gas. 

The model under consideration couples a Fokker-Planck-type equation on the sphere for the orientation distribution of the rods to the Navier-Stokes equations, which are now enhanced by additional stresses reflecting the orientation of the rods on the molecular level.  The coupled problem is 5- dimensional (three-dimensions in physical space and two degrees of freedom on the sphere) and it describes the interaction between the orientation of rod-like polymer molecules on the microscopic scale and the macroscopic properties of the fluid in which these molecules are contained. The macroscopic flow leads to a change of the orientation and, in the case of flexible particles, to a change in shape of the suspended microstructure. This process, in turn yields the production of a fluid stress.  The free-boundary problem is defined with the aid of a threshold for the pressure beyond which one has the incompressible Navier-Stokes equations for the fluid and below which one has a compressible model for the gas.

\subsection{Outline}
The outline of this article is as follows:
Section \ref{S1} presents  the main motivation for the upcoming investigation. Section \ref{S2} introduces modeling aspects of the problem: the physical setting, constitutive relations,  the free-boundary problem, and the statement of the problem which outlines the main objective of the work, namely the establishment of the global existence of the weak solutions to the free-boundary problem by  rigorously  showing that they can be obtained as the limit of weak solutions to the Doi problem for compressible fluids.
Additionally, the main results of the article as well as the notion of solutions to the macroscopic system are introduced.
Section \ref{S3} is dedicated to the construction of a  suitable approximate system and  the definition of its weak solution. Section \ref{S4} presents the existence of approximate solutions which relies on the derivation of suitable a priori estimates.  Section \ref{S5} presents the proof of the main theorem which relies on the establishment of the compactness of the sequence of solutions.

\section{Formulation of the main problem}\label{S2}
\subsection{Notations:} 
Before formulating the governing equation of our main problem we fix here some notations we are going to use in the paper.
\begin{itemize}
\item $L^{p}(0,T;X)$ denotes the Banach set of Bochner measurable functions $f$ from $(0,T)$ to $X$ endowed with either the norm $\Big(\int ^{T}_{0} \|g(\cdot, t)\|^{p}_{X}dt\Big)^{\frac{1}{p}}$ for $1\leq p<\infty$ or $\displaystyle\sup_{t>\infty} \|g(\cdot,t)\|_{X}$ for $p=\infty$. In particular, $f\in L^{r}(0, T; XY)$ denotes $\Big(\int ^{T}_{0} \big\|\big(\|f(t)\|_{Y_{\tau}}\big)\big\|^{p}_{X}dt\Big)^{\frac{1}{p}}$ or $\displaystyle\sup_{t>\infty} \big\|\big(\|f(t)\|_{Y_{\tau}}\big)\big\|_{X}$ for $p=\infty$. The notation
 $L^{p}_{t}L^{q}_{x}$ will abbreviate  the space $L^{p}(0,T);L^{q}(\Omega))$.
\item ${\mathcal M}((0,T) \times \Omega)$ is the space of bounded measures on $(0,T) \times \Omega$.
\item$A \lesssim B $ means there is a constant $C$ such that $A \leq C B $.
\item $\mathbf{1}_{X}$ is the indicator function which is 1 for $x\in X$ and 0 otherwise.
\item $C(T)$ is a function only depending on initial data and $T$, $C_{w}([0,T];X)$, is the space of continuos function from $(0,T)$ to $X$ endowed with the weak topology.
\item $\rightharpoonup$ and $\rightarrow$ denote weak limit and strong limit, respectively.
\item We denote by $\overline{x_{n}}$ the weak limit of a sequence $x_{n}$. 
\end{itemize}

\subsection{Governing equations}
We start by introducing the basic equations of motion for polymeric fluids. We recall  that a smooth motion of a body in continuum mechanics is described by a family of one-to-one mappings
\[
X(t, \cdot): \Omega \rightarrow \Omega, \quad t\in I.
\]
The curve $X(t,x)$ represents the trajectory of a particle occupying at time $t$ a spatial position $x$. A smooth motion $X$ is completely determined by a velocity field $u: I\times \Omega \rightarrow \mathbb{R}^{3}$ through
\[
\frac{\partial}{\partial t}X(t,x)=u(t, X(t,x)), \quad X(0,a)=a.
\]
Then, the conservation of mass can be formulated as follows:
\[
\frac{d}{dt}\int_{X(t,B)}\vr(t,x)dx=0, \quad B\subset \Omega.
\]
This equation is equivalent to
\[
\frac{d}{dt}\int_{B}\vr(t,x)dx +\int_{\partial B}\vr(t,x)[u(t,x)\cdot \hat{n}]dS=0,
\]
where $\hat{n}$ is the unit outer normal vector on $\partial \Omega$. If $\vr$ is smooth, one can use Green's theorem to deduce the following continuity equation:
\begin{equation} 
\label{eq:1.1}
\vr_{t} +\Div(\vr u)=0.
\end{equation}
We next obtain equation of motion by applying Newton's second law of motion.
\begin{equation} \label{eq:1.3}
\frac{d}{dt}\int_{B}\vr(t,x)u(u,x)dx+\int_{\partial B}\vr(t,x)[u(t,x)\cdot \hat{n}]dS=\int_{\partial B}\mathbb{T}(t,x)\hat{n}dS.
\end{equation}
By applying Green's lemma to (\ref{eq:1.3}), we finally have
\begin{equation} 
\label{eq:1.4}
(\vr u)_{t} +\Div(\vr u\otimes u)=\Div \mathbb{T}, \quad (\Div \mathbb{T})_{i} =\sum_{j=1}^{3}\frac{\partial \mathbb{T}_{ij}}{\partial x_{j}}.
\end{equation} 
The stress tensor $\mathbb{T}$ obeys Stokes' law:
\[
\mathbb{T}=\mathbb{S}-p\mathbb{I}.
\]
Let us determine $\mathbb{S}$ and $p$ in our model. $\mathbb{S}$ consists of two parts:
$$\mathbb{S}=\mathbb{S}_{1}+\mathbb{S}_{2}, $$
%
where $\mathbb{S}_{1}$ is the viscous stress tensor generated by the fluid 
$$
\mathbb{S}_{1}=\mu \big(\nabla u +(\nabla u)^{t}\big) +\lambda(\Div u)\mathbb{I},
$$
and  $\mathbb{S}_{2}$ is the macroscopic symmetric stress tensor derived from the orientation of the rods at the molecular level. 
The microscopic insertions at time $t$ and macroscopic place $x$ are described by the probability $f(t,x,\tau) d\tau.$
The suspension stress tensor $\mathbb{S}_2$ is given by an expansion
$$\mathbb{S}_2(x,t) = \sigma^{(1)}(x,t) + \sigma^{(2)}(x,t) + \sigma^{(3)}(x,t),$$
where
\[
\sigma^{(1)}(t,x)=\int_{S^{2}}(3\tau\otimes\tau-\mathbb{I}_{3\times3})f(t,x,\tau)d\tau, 
\]
\[
\sigma^{(2)}(t,x)  = -\sigma^{(2)}_{ij}(t,x) \mathbb{I}_{3 \times 3},\,\,\,  \mbox{with} \,\,\,\sigma^{(2)}_{ij}(t,x)= \int_{S^{2}} \gamma_{ij}^{(2)}(\tau) f(t,x,\tau) d \tau, 
\]
and 
\[
\sigma^{(3)}(t,x)  = -\sigma^{(3)}_{ij}(t,x) \mathbb{I}_{3 \times 3},\]
with
\[\sigma^{(3)}_{ij}(t,x)= \int_{S^{2}}  \int_{S^{2}} \gamma_{ij}^{(3)}(\tau_1, \tau_2)  f(t,x,\tau_1)  f(t,x, \tau_2) d \tau_1 d\tau_2.
\]
This, and more general expansions for $\mathbb{S}_2$ are encountered in the polymer literature (cf.\ Doi and Edwards  \cite{Doi1986}).
We refer the reader to the articles by Constantin et al \cite{Constantin2007}, \cite{Constantin2008}, where a general class of  stress tensors is presented in the context of incompressible fluids.

The structure coefficients in the expansion $\gamma_{ij}^{(2)}, \gamma_{ij}^{(3)}$ are in general smooth, time independent, $x$ independent, and do not depend on $f.$ Assuming for simplicity that 
$$\gamma_{ij}^{(2)} (\tau) = \gamma_{ij}^{(3)} (\tau_1, \tau_2) =1$$
and denoting 
 $$\eta(t,x)=\int_{S^{2}}f(t,x,\tau)d\tau$$
the suspension stress tensor $\mathbb{S}_2$  takes the form
\begin{equation} 
\label{eq:1.5}
\mathbb{S}_2(x,t) = \sigma^{(1)}(x,t) - \eta \mathbb{I}_{3 \times 3} -  \eta^2 \mathbb{I}_{3 \times 3}.
\end{equation}

In this setting, $f$ describes the time-dependent orientation distribution that a rod with a center mass at $x$ has an axis $\tau$ in the area element $d\tau$ and it is described by a compressible Fokker-Plank type equation,
\begin{equation} 
\label{eq:1.6}
f_{t}+\Div (f\vu)+\nabla_{\tau}\cdot(P_{\tau^{\perp}}\nabla u \tau f)-D_{\tau}\Delta_{\tau} f-D\Delta f=0,
\end{equation}
where $ P_{\tau^{\perp}}(\nabla_{x} u \tau)=\nabla_{x}u \tau-(\tau\cdot \nabla_{x}u \tau)\tau$ is the projection of $\nabla u \tau$ on the tangent space of $S^{2}$  at $\tau \in S^{2}$. With $\nabla_\tau$ and $\Delta_\tau$ we denote the gradient and the Laplace operator on the unit sphere, while $\nabla$ and $\Delta$ represent the gradient and the Laplacian operator in $\mathbb{R}^3$.

The second term $\nabla \cdot (uf)$ in \eqref{eq:1.6} describes the change of $f$ due to the displacement of the center of mass of the rods due to macroscopic advection. 
The term $\nabla_{\tau}\cdot(P_{\tau^{\perp}}\nabla u \tau f)$
is a drift-term on the sphere, which represents the shear-forces acting on the rods. 
The  term   $D_{\tau}\Delta_{\tau} f$  represents  the rotational  diffusion due to Brownian motion. This effect causes the rods to change their orientation spontaneously, whereas the term  $D\Delta f$  is the translational diffusion due to Brownian effects.
 
 By integrating (\ref{eq:1.6}) over $S^{2}$, we can obtain the equation of $\eta$:
\begin{equation} 
\label{eq:1.7}
\eta_{t}+\Div(\eta\vu)-D\Delta \eta=0.
\end{equation}
The pressure due to the fluid  $p_F$  is denoted by
\begin{equation} 
\label{eq:1.8}
p_F =\pi, 
\end{equation}
The pressure due to the presence of the dispersed particles denoted by $p_P$ is of the form
\begin{equation} 
\label{eq:1.9}
p_P = \eta + \eta^2.
\end{equation}
The overall pressure of the mixture is denoted by $P(\pi, \eta)$ and is given by
$$P(\pi, \eta) = p_F + p_P = \pi+ \eta + \eta^2.$$
By substituting (\ref{eq:1.5}) and (\ref{eq:1.8}) to (\ref{eq:1.4}), the equation of motion becomes
\begin{equation} 
\label{eq:1.9a}
\partial_t (\vr \vu)+\Div (\vr \vu \otimes \vu)-\mu \Delta \vu - \lambda \nabla \Div \vu+\nabla (\pi+\eta+ \eta^{2})=\Div \sigma^{(1)}.
\end{equation}
Finally, we have the following system of equations for the polymeric fluid in $(0,T)\times \Omega$:
\begin{equation} 
\partial_t \vr +\Div(\vr \vu)=0, \label{eq:1.10 a}
\end{equation}
\begin{equation}
\partial_t(\vr \vu) +\Div (\vr \vu \otimes \vu)-\Delta \vu - \nabla \Div \vu+\nabla P(\vr, \eta) =\Div \sigma^{(1)}, \label{eq:1.10 b}
\end{equation}
\begin{equation}
\partial_t f +\Div(f \vu) +\nabla_{\tau} \cdot(P_{\tau^{\perp}}(\nabla_{x} u \tau)f)- \Delta_{\tau}f- \Delta_{x}f=0, \label{eq:1.10 c}
\end{equation}
\begin{equation}
\label{eq:1.10 d}
\eta_{t}+\nabla\cdot(\eta\vu)-\Delta \eta=0.
\end{equation}
For the sake of simplicity, all the coefficients are normalized by 1. We consider (\ref{eq:1.10 a}), (\ref{eq:1.10 b}), \eqref{eq:1.10 c}, \eqref{eq:1.10 d} on a bounded domain with Dirichlet boundary conditions,
 $$u=0, \,\,  f=0,\, \mbox{and} \, \, \eta=0 \, \mbox{ on} \,\, \partial \Omega.$$
 \subsection{Statement of the problem}
In this article we are concerned with the free-boundary problem for the system \eqref{eq:1.10 a}- \eqref{eq:1.10 d} which is defined with the aid of a threshold for the pressure beyond which one has the incompressible Navier-Stokes equations for the fluid and below which one has a compressible model for the gas.

In sum, the free boundary problem {\bf(P)}  is governed by the following  equations in $(0,T)\times \Omega$:

\begin{equation} \label{2.14}
\partial_t \vr+\Div(\vr \vu)=0
 \end{equation}
\[
\quad  0 \le \vr \le 1 \nonumber
\]
\begin{equation}\label{2.15}
\partial_t (\vr  \vu) +\Div(\vr \vu\otimes \vu)+ \nabla(\pi + \eta + \eta^2) = \Div \mathbb{S} + \Div \sigma 
\end{equation}
\begin{equation}\label{2.16}
 \partial_t f+\Div(f \vu) +\nabla_{\tau} \cdot(P_{\tau^{\perp}}(\nabla_{x} \vu \tau)f)- \Delta_{\tau}f- \Delta_{x}f=0 
\end{equation}
\begin{equation}\label{2.17}
\partial_t \eta + \Div(\eta \vu)  - \Delta \eta = 0
\end{equation}
and the  free boundary conditions
\begin{equation}\label{2.18}
\Div \vu = 0 \quad \mbox{a.e} \,\, \mbox{on} \,\, \{ \vr = 1\}
\end{equation}
\begin{equation}\label{2.19}
\pi \ge  0  \quad \mbox{a.e} \,\, \mbox{in} \,\, \{ \vr = 1\}
\end{equation}
\begin{equation}\label{2.20}
\pi =  0  \quad \mbox{a.e} \,\, \mbox{in} \,\, \{ \vr < 1\}
\end{equation}
The unknowns here are the density $\vr,$ the velocity vector field $\vu,$ the pressure $\pi, $ which is Lagrange multiplier associated with the incompressibility constraint \eqref{2.18}  $\Div \vu = 0$ a.e. in $\{\vr = 1\},$ the orientation distribution $f$ and the particle density $\eta$. Note that $\pi$ is apparent only in the congested regions $\{\vr =1\}.$ In fact conditions \eqref{2.19}, \eqref{2.20},  can be rewritten, in an equivalent way, as one constraint
\begin{equation}\nonumber
\vr \pi = \pi\ge 0.
\end{equation}
The physical properties of the mixture are reflected through the following constitutive relations.
\medskip 

\begin{center}
{\bf Constitutive relations}
\end{center}
\vspace{-0.1in}
\begin{equation}\nonumber
 P(\pi, \eta) = p_F + p_P = \pi + \eta + \eta^2.
\end{equation}
\begin{equation}\nonumber
\mathbb{S} = \mu (\nabla \vu + \nabla \vu^{T}) + \xi \Div \vu \mathbb{I}
\end{equation}
\begin{equation}\nonumber
\Div \mathbb{S}=\mu \Delta \vu - \xi \nabla \Div \vu
\end{equation}
\begin{equation}\nonumber
\sigma(t,x) = \sigma^{(1)}(t,x)=\int_{S^{2}}(3\tau\otimes\tau-\mathbb{I}_{3\times3})f(t,x,\tau)d\tau,
\end{equation}
where,  for the sake of simplicity, all the coefficients will be normalized by 1. 
\medskip

\begin{center}
{\bf Boundary conditions}
\end{center}

We consider the problem {\bf(P)}  on a bounded domain with Dirichlet boundary conditions.

 \begin{equation}\label{2.22}
 \vu=0, \,\,  f=0,\, \mbox{and} \, \, \eta=0, \, \mbox{ on} \,\  \partial \Omega.
 \end{equation}
 
 \medskip

\begin{center}
{\bf Initial data}
\end{center}
The system must be complemented with initial conditions, namely
\begin{equation}\label{2.23}
 \vr|_{t=0} = \vr_{0},\quad \vr \vu|_{t=0} = m_{0}, \quad \eta|_{t=0}=\eta_{0},\quad f|_{t=0}=f_{0}
\end{equation}
where
\begin{eqnarray}
& & 0 \le \vr_{0} \le 1, \quad  \vr_{0} \in L^1(\Omega)\nonumber\\
& & m_{0} \in L^2 (\Omega), \quad m_{0} =0 \,\, \mbox{a.e.  on }\,\, \{\vr_{0} = 0\}\nonumber\\
& & \vr_{0} \not \equiv 0,\,\,\, \int_{\Omega}\!\!\!\!\!\!\!-\vr_{0}=M<1\label{2.29}\nonumber \\
&& \vr_{0} |\vu_{0}|^{2}  \in L^2(\Omega),\quad \vu^0 = \frac{m_{0}}{\vr_{0}}\,\, \mbox{on} \,\, \{\vr_{0} >0\}\nonumber \\
&& \eta_{0}\in L^{2}(\Omega)\nonumber\\
&& f_{0}\in L^{1}(\Omega\times S^{2})\nonumber
\end{eqnarray}

The goal of   this paper is to prove the  existence of weak solutions to the free-boundary problem \eqref{2.14}-\eqref{2.20}, so we introduce the notion of weak solutions we are going to use through the paper.
 
 \begin{definition}
\label{def1}
{\bf [Weak solution of the problem (P)]}\\
A vector $(\vr, \vu,\pi, f, \eta)$ is called a weak solution to \eqref{2.14}-\eqref{2.20} with boundary data \eqref{2.22} and  initial data \eqref{2.23} if the equations
\begin{equation} \nonumber
\partial_t \vr+\Div(\vr \vu)=0
 \end{equation}
\begin{equation}\nonumber 
\partial_t (\vr  \vu) +\Div(\vr \vu\otimes \vu)+ \nabla(\pi + \eta + \eta^2) = \Div \mathbb{S} + \Div \sigma
\end{equation}
\begin{equation}\nonumber
 \partial_t f+\Div(f \vu) +\nabla_{\tau} \cdot(P_{\tau^{\perp}}(\nabla_{x} \vu \tau)f)- \Delta_{\tau}f- \Delta_{x}f=0 
\end{equation}
\begin{equation}\nonumber
\partial_t \eta + \Div(\eta \vu)  - \Delta \eta = 0
\end{equation}
are satisfied in the sense of distributions, the divergence free condition $\Div \vu = 0$ is satisfied a.e. in   $\{ \vr = 1\},$
the constrained $0\le \vr \le 1$ is satisfied a.e. in $(0,T) \times \Omega$ and the following regularity properties hold
\[\vr \in C([0,T]; L^p(\Omega)), \,\, 1\le p < \infty,\]
\[\vu \in L^2(0,T;(W_0^{1,2}(\Omega))),\,\, \vr|{\bf u}|^2 \in L^{\infty}(0, T;L^1(\Omega)),\]
\[\pi \in {\mathcal M}((0,T) \times \Omega)\]
\[\eta \in L^{\infty}(0,T; L^{2}(\Omega))\cap L^{2}(0,T; \dot{H}^{1}(\Omega)), \quad f\ln f \in L^{\infty}(0,T; L^{1}(\Omega\times S^{2}))\]
\[ \nabla _{\tau}\sqrt{f} \in L^{2}((0,T)\times\Omega\times S^{2} ), \quad \nabla\sqrt{f} \in L^{2}((0,T)\times\Omega\times S^{2}).\]
\end{definition}
Moreover $\pi $ is so regular that the condition 
\[\pi(\vr-1)=0, \]
is satisfied in the sense of distribution.
The objective of this work is to prove the existence of weak solutions to the free-boundary problem \eqref{2.14}-\eqref{2.20} by showing rigorously that they   can be obtained as a limit of $(\vr_n, \vu_n, f_n, \eta_n)$ - the weak solutions to the Doi model for compressible fluids 
\begin{equation} 
\partial_t \vr_n +\Div(\vr_n \vu_n)=0, \label{2.32}
\end{equation}
\begin{equation}
\partial_t(\vr_n \vu_n) +\Div (\vr_n \vu_n \otimes \vu_n)-\Delta \vu_n-\nabla \Div \vu_n+\nabla (\pi_n+ \eta_n+ \eta_n^2) =\Div\sigma_n, \label{2.33}
\end{equation}
\begin{equation}
\partial_t f_n +\Div(f_n \vu_n) +\nabla_{\tau} \cdot(P_{\tau^{\perp}}(\nabla_{x} \vu_n \tau)f)- \Delta_{\tau}f_n- \Delta_{x}f_n=0, \label{2.34}
\end{equation}
\begin{equation}
\label{2.35}
\partial_t \eta_n+\nabla\cdot(\eta_n\vu_n)-\Delta \eta_n=0,
\end{equation}
where
$$\pi_{n}=(\vr_{n})^{\gamma_{n}},\  \gamma_{n}\to \infty,\ \text{as $n\to \infty$}.$$

Finally, we want to point out that the solution we are going to obtain satisfies the following energy inequality
\begin{equation}\label{energy}
 \begin{split}
 &\int_{\Omega} \Big(\frac{\rho|\vu|^{2}}{2}+\eta^{2} +\psi \Big)(t)dx + 4\int^{t}_{0}\int_{\Omega}\int_{S^{2}} |\nabla_{\tau}\sqrt{f}|^{2}d\tau dxdt \\
 &+4 \int^{t}_{0} \int_{\Omega} \int_{S^{2}} |\nabla \sqrt{f}|^{2}d\tau dxdt +\int^{t}_{0} \int_{\Omega}\Big(|\nabla \vu|^{2} + |\Div \vu|^{2} +2|\nabla \eta|^{2}\Big)dxdt \\
 & \leq \int_{\Omega} \Big(\frac{\rho_{0}|u_{0}|^{2}}{2}+\eta^{2}_{0} +\psi_{0} \Big)dx,
 \end{split}
\end{equation}
where
$$\psi(t,x)=\int_{S^{2}} (f\ln f)(t,x,\tau) d\tau.$$

\subsection{Main results}
Now we are ready to state the main existence results for our problem.
\begin{theorem}
Assume that the boundary conditions \eqref{2.22} and the initial conditions \eqref{2.23} are satisfied. Then,  there exists a weak solution (in the sense of Definition \ref{def1}) of the problem \eqref{2.14}-\eqref{2.20}. 
\label{MT}
\end{theorem}
The main Theorem \ref{MT} will be obtained  as  a consequence of the following result.

\begin{theorem} 
Let $n \in \N$ be fixed, then there exists a global weak solution  $(\vr_n, \vu_n, \pi, f_n, \eta_n)$ to 
\eqref{2.32}-\eqref{2.35} in the sense of the Definition \ref{def2}, such that, as $n \to \infty$
\begin{equation}
\label{2.1.1}
(\vr_{n}-1)_{+}\rightarrow 0\qquad \text{in $L^{\infty}(0,T;L^{p})$, for any $1\leq p\leq 0.$}
\end{equation}
Moreover,
\begin{equation}
(\vr_{n})^{\gamma_{n}} \qquad \text{is bounded in $L^{1}$, for $n$ such that $\gamma_{n}\geq 3$,}
\end{equation}
and up to a subsequence there exists $\pi\in \mathcal{M}((0,T)\times \Omega)$ such that
\begin{equation}
(\vr_{n})^{\gamma_{n}} \rightharpoonup \pi, \qquad \text{as $n\to \infty$}.
\end{equation}
If in addition $\vr_{n0}\to\vr_{0}$ in $L^{1}$, then the following convergence holds:
 \[
 \vr_n \rightharpoonup \vr \,\, \mbox{weakly in} \,\, L^p((0,T) \times \Omega) \,\, 1\le p < +\infty,
 \]
 \[
\vr_n \vu_n \rightharpoonup \vr\vu  \,\, \mbox{weakly in} \,\, L^p((0,T; L^{r}(\Omega)),\  1\le p < +\infty,\ 1\leq r<2,
 \]
 \[
\vr_n \vu_n\otimes \vu_n \rightharpoonup \vr\vu\otimes\vu \,\, \mbox{weakly in} \,\, L^p((0,T; L^{1}(\Omega)),1\le p < +\infty,
 \]
 \[
 f_n \rightharpoonup f \,\, \mbox{weakly in} \,\, L^2((0,T; L^{6/5}(\Omega\times S^{2})),
 \]
 \[
 \eta_n \rightarrow \eta  \,\, \mbox{strongly in} \,\, L^{2}(0,T;L^{2}(\Omega)),
 \]
 $0\le \vr \le 1$
and 
 $(\vr,  \vu, \pi, f, \eta)$ is a  weak solution to the problem \eqref{2.14}-\eqref{2.20} in the sense of Definition \ref{def1}. 
 \label{MT2}
 \end{theorem}
 
The rest of the paper is devoted to the proof of the Theorems \ref{MT} and  \ref{MT2}.

\section{Approximating problem}\label{S3}
We describe, now,    the approximating scheme we are going to use. 

Let be $\gamma_{n}$ a sequence of real numbers such that $\gamma_{n}>\frac{3}{2}$, for any $n\in \N$ and $\gamma_{n}\to \infty$ as 
$n\to \infty$, we define  $\{\vr_{n}, \vu_{n}, f_{n}, \eta_{n}\}$ as solutions of the following system

\begin{equation} \label{4.1}
\partial_t \vr_n+\Div(\vr_n \vu_n)=0, \quad  \vr_n \ge 0
 \end{equation}
\begin{equation}\label{4.2}
\partial_t (\vr_n  \vu_n) +\Div(\vr_n \vu_n\otimes \vu_n) - \Delta \vu_n - \nabla \Div \vu_n + \nabla (\pi_{n} + \eta_n + \eta_n^2)   =  \Div\sigma_n 
\end{equation}
\begin{equation}\label{4.3}
 \partial_t f_n+\Div(f_n \vu_n) +\nabla_{\tau} \cdot(P_{\tau^{\perp}}(\nabla_{x} \vu_n \tau)f_n)- \Delta_{\tau}f_n- \Delta_{x}f_n=0 
\end{equation}
\begin{equation}\label{4.4}
\partial_t \eta_n + \Div(\eta_n \vu_n)  - \Delta \eta_n = 0,
\end{equation}
where
$$\pi_{n}=(\vr_n)^{\gamma_n}$$
and 
\begin{equation}\label{4.5}
\sigma_n(t,x) = \int_{S^{2}}(3\tau\otimes\tau-\mathbb{I}_{3\times3})f_n(t,x,\tau)d\tau
\end{equation}
The approximating system must be complemented with boundary and initial data as follows.
\begin{center}
{\bf Boundary data}
\end{center}
\begin{equation}
 \label{b1}
 \vu_{n}=0, \,\,  f_{n}=0,\, \mbox{and} \, \, \eta_{n}=0, \, \mbox{ on} \,\  \partial \Omega.
 \end{equation}
\vspace{0.1in}
\begin{center}
{\bf Initial data}
\end{center}
\begin{equation}
\vr_n|_{t=0} = \vr_{n_0}, \quad \vr_n \vu_n|_{t=0} = m_{n_0}, \quad \eta_{n}|_{t=0}=\eta_{n_0},\quad f_{n}|_{t=0}=f_{n_0}
\label{i1}
\end{equation}
where
\begin{eqnarray}
& & 0 \le \vr_{n_0} \quad \mbox{a.e}, \quad  \vr_{n_{0}} \in L^1(\Omega) \cap L^{\gamma_n}(\Omega), \nonumber\\
& & \int (\vr_{n_0})^{\gamma_n} dx \le c \gamma_n \,\, \mbox{for some}\,\, c, \label{i5}\\
& & \quad \quad \quad  m_{n_0} \in L^{\frac{2 \gamma_n}{\gamma_n + 1}} (\Omega), \nonumber \\
&& \vr_{n_0} |\vu_{n_0}|^2  \,\, \mbox{is bounded in}\,\,  L^1(\Omega), \nonumber\\
&& \quad \quad \quad \vu_{n_0} = \frac{m_{n_0}}{\vr_{n_0}}\,\, \mbox{on} \,\, \{\vr_{n_0} >0\}, \nonumber\\
&&  \quad \quad \quad \vu_{n_0} = 0\,\, \mbox{on} \,\, \{\vr_{n_0} = 0\}, \nonumber\\
&&  \quad \quad \quad  f_{n_0}  \in L^1(\Omega \times S^{2}), \nonumber\\
&&  \quad \quad \quad \eta_{n_0}  \in L^2(\Omega \times S^{2}).\nonumber
\end{eqnarray}
Furthermore we assume that
\begin{equation}
M_{n}=\int_{\Omega}\!\!\!\!\!\!\!-\vr_{n_0},\quad  0<M_{n}<M<1,\quad M_{n}\to M.
\label{i4}
\end{equation}
\begin{equation}\nonumber
\vr_{n_0}\vu_{n}\rightharpoonup m_{0}\quad \text{weakly in $L^{2}(\Omega)$,}
\end{equation}
\begin{equation}
\vr_{n_0}\rightharpoonup \vr_{0}\quad \text{weakly in $L^{1}(\Omega)$.}
\label{i3}
\end{equation}
\subsection{Definition of weak solution of the approximate system} \label{sec:2.2}
For any fixed $\gamma_{n}>3/2$ we now define the notion of  weak solution of the system \eqref{4.1}-\eqref{4.4}, with initial data \eqref{i1}  and boundary data \eqref{b1}. 
\begin{definition} \label{def2}
For any fixed $\gamma_{n}>3/2$, we say $\{\vr_n, \vu_n, f_n, \eta_n, \sigma_n\}$ is a weak solution of the system \eqref{4.1}-\eqref{4.4} if
\begin{enumerate}[]
\item (i) 
$$\vr_{n} \in L^{\infty}(0,T; L^{\gamma_{n}}(\Omega)), \quad \nabla \vu_n \in L^{2}(0;T; L^{2}(\Omega)),$$
$$\vr_n|\vu_n|^{2}\in L^{\infty}(0,T; L^{1}(\Omega)), \quad \vr_n \vu_n \in C_{w}([0,T];L^{\frac{2\gamma_{n}}{\gamma_{n}+1}}(\Omega)),$$
$$\eta_n \in L^{\infty}(0,T; L^{2}(\Omega))\cap L^{2}(0,T; \dot{H}^{1}(\Omega)), \quad f_n\ln f_n \in L^{\infty}(0,T; L^{1}(\Omega\times S^{2}))$$
$$\nabla _{\tau}\sqrt{f_n} \in L^{2}((0,T)\times\Omega\times S^{2}), \quad \nabla\sqrt{f_n} \in L^{2}((0,T)\times\Omega\times S^{2}),$$
\item (ii) \eqref{4.1} holds in the sense of renormalized solutions, i.e.,
\begin{equation} \label{eq:2.10}
\partial_{t}(b(\vr_n)) +\Div (b(\vr_n)\vu_n)+\left(b'(\vr_n)\vr_n-b(\vr_n)\right)\Div \vu_n=0
\end{equation}
holds in the sense of distributions for any $b\in C^{1}$ such that $|b^{'}(z)z|+|b(z)| \leq C$ for all $z\in \mathbb{R}$,
\item (iii) \eqref{4.2}, \eqref{4.3}, \eqref{4.4}, and \eqref{4.5} hold in the sense of distributions,
\item (iv)  the following energy inequality is satisfied :
\begin{equation}\label{4.19}
 \begin{split}
 &\int_{\Omega} \Big[\frac{\vr_n |\vu_n|^{2}}{2}+\frac{\vr_n^{\gamma_{n}}}{\gamma_{n}-1}+\eta_n^{2} +\psi_n \Big](t) 
 dx + 4 \int_0^t \int_{\Omega}\int_{S^{2}} |\nabla_{\tau}\sqrt{f_n}|^{2}d\tau dx dt + \\
 &4 \int_0^T \int_{\Omega} \int_{S^{2}} |\nabla \sqrt{f_n}|^{2}d\tau dx dt + \int_0^t \int_{\Omega}\Big[|\nabla \vu_n|^{2} + |\Div \vu_n|^{2} +2|\nabla \eta_n|^{2}\Big]dx dt \leq\\
 & \int_{\Omega} \Big[\frac{\vr_{n_0} |\vu_{n_0}|^{2}}{2}+\frac{\vr_{n_0}^{\gamma_{n}}}{\gamma_{n}-1}+\eta_{n_0}^{2} +\psi_{n_0} \Big] dx=E_{n_{0}}, 
 \end{split}
\end{equation}
where
\[
\psi_n(t,x)=\int_{S^{2}} (f_n\ln f_n)(t,x,\tau) d\tau.
\]
\end{enumerate}
\end{definition}

\section{Existence of approximate solutions}\label{S4}
 For any fixed $n\in \N$, the existence of weak solutions for the system \eqref{4.1}-\eqref{4.4} has been proved by Bae and Trivisa in \cite{BT2012} (we refer the reader to \cite{BT2013} for the treatment of the Doi model for incompressible-polymeric fluids) , we can summarize their existence result as follows.
\begin{theorem} \label{thm:2.2}
Let $\gamma_{n}>\frac{3}{2}$ and $\Omega$ be a $C^{1}$ bounded domain. Assume that the initial data $\{\vr_{n_0}, \vu_{n_0}, f_{n_0}, \eta_{n_0}\}$ satisfy \eqref{i4}-\eqref{i3} and the boundary conditions \eqref{b1} hold. Then, there exists a weak solution  (in the sense of Definition \ref{def2}) $\{\vr_n, \vu_n, f_n, \eta_n, \sigma_n\}$ of the system \eqref{4.1}-\eqref{4.4} satisfying \eqref{i1} at $t=0$. 
\end{theorem}
By following the same line of arguments of \cite{BT2012} we recall in the next section the main compactness properties of the approximate solution $\{\vr_n, \vu_n, f_n, \eta_n, \sigma_n\}$.
\subsection{Energy estimates of the approximating system}
Besides the bounds mentioned in (i) of Definition \ref{def2} we can collect some further estimates satisfied by the solutions $\{\vr_n, \vu_n, f_n, \eta_n, \sigma_n\}$.
By the energy inequality \eqref{4.19} and the Sobolev embedding $\dot{H}^{1}\subset L^{6}$, we can estimate $\sqrt{f_{n}}$ as
\[
\sqrt{f_{n}} \in L^{2}\big(0,T; L^{2}(\Omega)L^{6}(S^{2}) \cap L^{6}(\Omega)L^{2}(S^{2})\big).
\]
This implies that
\begin{equation} \label{eq:2.7}
f_{n}\in L^{1}\big(0,T; L^{1}(\Omega)L^{3}(S^{2}) \cap L^{3}(\Omega)L^{1}(S^{2})\big) \subset L^{1}(0,T; L^{2}(\Omega\times S^{2})).
\end{equation}
We finally estimate $\sigma_{n}$. Since $|\sigma_{n}(t,x)| \leq 3 \displaystyle\int_{S^{2}}f_{n}(t,x,\tau) d\tau=3\eta_{n}(t,x)$,
\begin{equation} \nonumber
\sigma_{n} \in L^{1}(0,T; L^{3}(\Omega)) \cap L^{\infty}(0,T; L^{2}(\Omega))
\end{equation}
where the first space is derived by integrating $f_{n}$ over $S^{2}$ using (\ref{eq:2.7}) and the second bound is from $\eta_{n} \in L^{\infty}(0,T; L^{2}(\Omega))$. We next estimate the derivative of $\sigma$. By using the entropy dissipation,
\begin{equation*}
 \begin{split}
 |\nabla \sigma_{n}(t,x)| &\leq 3 \int_{S^{2}}|\nabla f_{n}(t,x,\tau)|d\tau \lesssim \left[\int_{S^{2}} |\nabla \sqrt{f_{n}}|^{2}d\tau \right]^{\frac{1}{2}} \left[\int_{S^{2}}(\sqrt{f_{n}})^{2}d\tau \right]^{\frac{1}{2}}\\
 &=\left[\int_{S^{2}} |\nabla\sqrt{f_{n}}|^{2}d\tau \right]^{\frac{1}{2}} (\eta_{n})^{\frac{1}{2}}.
 \end{split}
\end{equation*}
Since
\[
(\eta_{n})^{\frac{1}{2}} \in L^{\infty}(0,T; L^{4}(\Omega)) \cap L^{2}(0,T; L^{6}(\Omega)), \quad \left[\int_{S^{2}} |\nabla\sqrt{f_{n}}|^{2}d\tau \right]^{\frac{1}{2}} \in L^{2}(0,T; L^{2}(\Omega)),
\]
we have
\begin{equation}\nonumber
\nabla \sigma_{n} \in L^{1}(0,T; L^{\frac{3}{2}}(\Omega)) \cap L^{2}(0,T; L^{\frac{4}{3}}(\Omega)).
\end{equation}
Moreover, as we will see in the Section \ref{S4} we will be able to show to following uniform bound in $n\in \N$ for $\vr_{n}$
\begin{equation}\nonumber
\vr_{n}\in  L^{\infty}(0,T; L^{1}\cap L^{p}(\Omega)),\qquad 1\leq p<+\infty.
\end{equation}
 Extracting a subsequence, using the same notation, $\{\vr_{n}, \vu_{n}, f_{n}, \eta_{n}, \sigma_{n}\}_{n \ge 1}$, we have various limit functions such as
\begin{equation} 
 \begin{split}
 & \vr_{n} \rightharpoonup \vr \hspace{0.2cm} \text{in} \hspace{0.2cm} L^{\infty}(0,T;L^{p}(\Omega)), \quad \vr \in L^{\infty}(0,T; L^{1}\cap L^{p}(\Omega)),  \hspace{0.2cm} 1\leq p<+\infty,\\
 & \sqrt{\vr_{n}} \vu_{n} \rightharpoonup v \hspace{0.2cm} \text{in} \hspace{0.2cm} L^{2}(0,T; L^{2}(\Omega)), \quad v \in L^{\infty}(0,T; L^{2}(\Omega)),\\
 & \vu_{n} \rightharpoonup \vu \hspace{0.2cm} \text{in} \hspace{0.2cm} L^{2}(0,T; H^{1}(\Omega)),\\
 & \vr_{n} \vu_{n} \rightharpoonup m \hspace{0.2cm} \text{in} \hspace{0.2cm} L^{\frac{2p}{p+1}}(\Omega\times(0,T)), \quad m \in L^{\infty}(0,T; L^{\frac{2p}{p+1}}(\Omega)),  \hspace{0.2cm} 1\leq p<+\infty, \\
 & \vr_{n}\vu_{n_i}\vu_{n_j} \rightharpoonup e_{ij} \hspace{0.2cm} \text{in the sense of measures}, \\
 &\hspace{3cm} e_{ij} \hspace{0.2cm} \text{is a bounded measure}, \\
 & f_{n} \rightharpoonup f \hspace{0.2cm} \text{in} \hspace{0.2cm} L^{2}(0,T; L^{\frac{6}{5}}(\Omega\times S^{2})),\\
 & \eta_{n} \rightharpoonup \eta \hspace{0.2cm} \text{in} \hspace{0.2cm} L^{2}(0,T; H^{1}(\Omega)),\quad \eta \in L^{\infty}(0,T; L^{2}(\Omega)) \cap L^{2}(0,T; H^{1}(\Omega)),\\
 & \sigma_{n} \rightharpoonup \sigma \hspace{0.2cm} \text{in} \hspace{0.2cm} L^{2}(0,T; L^{2}(\Omega)), \quad \sigma \in L^{\infty}(0,T; L^{2}(\Omega)) \cap L^{1}(0,T; L^{3}(\Omega)),\\
 &\nabla \sigma_{n} \rightharpoonup \nabla \sigma \hspace{0.1cm} \text{in} \hspace{0.1cm} L^{2}(0,T; L^{\frac{4}{3}}(\Omega)), \quad \nabla \sigma \in L^{1}(0,T; L^{\frac{3}{2}}(\Omega)) \cap L^{2}(0,T; L^{\frac{4}{3}}(\Omega)).
 \end{split}
  \label{eq:2.16}
\end{equation}

Finally  we state here the following compactness results (for the proof we refer to \cite{BT2012}, Proposition 2.1)
\begin{proposition}[Compactness] \label{prop:2.3}
The limit functions in (\ref{eq:2.16}) satisfy the following statements.
\begin{enumerate}[]
\item (i) $v=\sqrt{\vr}\vu$, $m=\vr \vu$, $e_{ij}=\vr \vu_{i} \vu_{j}$.
\item (ii) $\eta_{n}$ converges strongly to $\eta$ in $L^{2}(\Omega\times (0,T))$, and $\sigma_{n}$ converges strongly to $\sigma$ in $L^{2}(\Omega \times (0,T))$.
\item (iii) $\vr_{n}(\eta_{n})^{2}$ converges to $\vr \eta^{2}$ in the sense of distributions.
\item (iv) $\vr$ and $\vu$ solve (\ref{eq:1.10 a}) in the sense of renormalized solutions.
\item (v) If in addition we assume that $\vr_{n0}$ converges to $\vr_{0}$ in $L^{1}(\Omega)$,
\begin{equation}\nonumber
\vr_{n} \rightarrow \vr \hspace{0.2cm} \text{in} \hspace{0.2cm} L^{1}(\Omega\times(0,T)) \cap C([0,T]; L^{p}(\Omega)) \hspace{0.2cm} \text{for all} \hspace{0.2cm} 1\leq p<+\infty.
\end{equation}
\item (vi) Finally, we have the following strong convergence:
\begin{equation} \nonumber
 \begin{split}
 & \vr_{n}\vu_{n} \rightarrow \vr \vu \hspace{0.2cm} \text{in} \hspace{0.2cm} L^{p}(0,T; L^{r}(\Omega)) \hspace{0.2cm} \text{for all} \hspace{0.2cm} 1\leq p<\infty, \quad 1\leq r<2, \\
 & \vu_{n} \rightarrow \vu \hspace{0.2cm} \text{in} \hspace{0.2cm} L^{q}(\Omega\times (0,T))\cap\{\vr_{n}>0\} \hspace{0.2cm} \text{for all} \hspace{0.2cm} 1\leq q<2, \\
 & \vu_{n} \rightarrow \vu \hspace{0.2cm} \text{in} \hspace{0.2cm} L^{2}(\Omega\times (0,T))\cap\{\vr_{n} \ge \delta\} \hspace{0.2cm} \text{for all} \hspace{0.2cm} \delta>0, \\
 & \vr_{n}\vu_{ni}\vu_{nj} \rightarrow \vr \vu_{i}\vu_{j} \hspace{0.2cm} \text{in} \hspace{0.2cm} L^{p}(0,T; L^{1}(\Omega)) \hspace{0.2cm} \text{for all} \hspace{0.2cm} 1\leq p<\infty.
 \end{split}
\end{equation}
\end{enumerate}
\end{proposition}

\section{Proof of the main theorem \ref{MT}}\label{S5}
This section is devoted to the proof of the Main Theorem \ref{MT}. We start with the proof Theorem \ref{MT2}, since, as we will see later on the Theorem \ref{MT} is a consequence of it.
\subsection{Proof of the  Theorem  \ref{MT2}}

 For simplicity we divide the proof in different steps.\\[0.5cm]
{\bf Step 1: Convergence of  $\mathbf{(\vr_{n}-1)_{+}}$ to $0$.}\\
By combining  the energy inequality \eqref{4.19} with \eqref{i5} we obtain
\begin{equation}
\label{e1}
\int_{\Omega}(\vr^{\gamma_{n}})_{n} dx\leq (\gamma_{n}-1)E_{n_{0}}+\int_{\Omega} (\vr_{n_0})^{\gamma_n} dx \leq (\gamma_{n}-1)E_{n_{0}}+c \gamma_n
\leq c\gamma_{n}.
\end{equation}
Since $\gamma_{n}\to \infty$ there exists  $n\in \N$ such that $\gamma_{n}>p$,  $1<p<+\infty$, then by H\"older inequality we get
\begin{equation}\nonumber
\|\vr_{n}\|_{L^{\infty}_{t}L^{p}_{x}}\leq \|\vr_{n}\|_{L^{\infty}_{t}L^{1}_{x}}^{\theta_{n}} \|\vr_{n}\|_{L^{\infty}_{t}L^{1}_{x}}^{1-\gamma_{n}}\leq
M_{n}^{\theta_{n}}(c\gamma_{n})^{\frac{1-\theta_{n}}{\gamma_{n}}},
\end{equation}
where $M_{n}$ is defined in \eqref{i4} and $\displaystyle{\frac{1}{p}=\theta_{n}+\frac{1-\theta_{n}}{\gamma_{n}}}$. As $n\to \infty$ we have that $\displaystyle{\theta_{n}\to \frac{1}{p}}$ and
\begin{equation}\nonumber
\|\vr_{n}\|_{L^{\infty}_{t}L^{p}_{x}}\leq \liminf_{n\to\infty}\|\vr_{n}\|_{L^{\infty}_{t}L^{p}_{x}}\leq M^{1/p}.
\end{equation}
Let us define the function $\phi_{n}$ as follows
$$\phi_{n}=(\vr_{n}-1)_{+},$$
by using again the energy inequality \eqref{4.19} we can compute
\begin{equation}
\label{e4}
\int_{\Omega}(1+\phi_{n})^{\gamma_{n}}\mathbf{1}_{\{\phi_{n}>0\}}dx \leq \int_{\Omega}\vr^{\gamma_{n}} dx\leq c\gamma_{n}.
\end{equation}
We recall the following inequality
\begin{equation}\nonumber
(1+x)^{k}\geq 1+c_{p}k^{p}x^{p}, \quad p>1, k \ \text{large}
\end{equation}
that holds for  any nonnegative $x$, and we apply it with $k=\gamma_{n}$, $x=\phi_{n}$ to the right hand side of \eqref{e4}, so we get
\begin{equation}\nonumber
c_{p}\gamma_{n}^{p}\int_{\Omega}\phi_{n}^{p}dx\leq |\Omega|+c_{p}\gamma_{n}^{p}\int_{\Omega}\phi_{n}^{p}dx\leq\int_{\Omega}(1+\phi_{n})^{\gamma}\mathbf{1}_{\{\phi_{n}>0\}}dx\leq c\gamma_{n}
\end{equation}
Hence we have
\begin{equation}\nonumber
\int_{\Omega}\phi_{n}^{p}dx\leq \frac{c}{c_{p}\gamma_{n}^{p-1}},
\end{equation}
and, as $n\to \infty$ we obtain
\begin{equation}\nonumber
(\vr_{n}-1)_{+}\rightarrow 0 \qquad \text{in $L^{\infty}(0,T;L^{p}(\Omega))$,  $1\leq p<+\infty$.}
\end{equation}\\[0.5cm]
{\bf Step 2: $\mathbf{L^{1}}$ uniform bound of $\mathbf{(\vr_{n})^{\gamma_{n}}}$.}\\
Assume that we know 
\begin{equation}
(\vr_{n})^{\gamma_{n}+1} \quad \text{is uniformly bounded in $L^{1}(0,T;L^{1}(\Omega))$},
\label{e8}
\end{equation}
then we have
\begin{equation}
\begin{split}
\int_{0}^{T}\!\!\int_{\Omega}(\vr_{n})^{\gamma_{n}}dxdt&=\int_{0}^{T}\!\!\left(\int_{\Omega\cap\{\vr_{n}>1\}}(\vr_{n})^{\gamma_{n}}dx+\int_{\Omega\cap\{\vr_{n}\leq1\}}(\vr_{n})^{\gamma_{n}}dx\right)dt\\
&\leq \int_{0}^{T}\!\!\left(\int_{\Omega}\left((\vr_{n})^{\gamma_{n}+1} +\vr_{n}\right)dx \right)dt.
\end{split}
\label{e9}
\end{equation}
By using \eqref{e8} and the fact that $\vr_{n}\in L^{\infty}(0,T;L^{1}(\Omega))$, from \eqref{e9} it follows the uniform $L^{1}$ bound for $(\vr_{n})^{\gamma_{n}}$.

In order to complete this step we have only to prove \eqref{e9}. We recall  that for $\vr_{n}$ we don't have $L^{\infty}$ bounds, but on the other hand, 
because of \eqref{e1} there exists a constant $\tilde{c}$ such that for any $n\in \N$ the following estimate holds
\begin{equation}
\|\vr_{n}\|_{L^{\infty}_{t}L^{\gamma_{n}}_{x}}\leq \tilde{c},
\label{e10}
\end{equation}
where $\displaystyle{\tilde{c}=\sup_{\gamma>0}(c\gamma)^{1/\gamma}}$.

We  define, now, the operator $\mathcal{B}$ as the inverse of the divergence operator. We denote the solution $v$ of
\[
\Div v=g \hspace{0.2cm} \text{in} \hspace{0.2cm} \Omega, \quad v=0 \hspace{0.2cm} \text{on} \hspace{0.2cm} \partial \Omega.
\]
by $v=\mathcal{B}g$. The operator $\mathcal{B}=( \mathcal{B}_{1}, \mathcal{B}_{2}, \mathcal{B}_{3})$ is the inverse of the divergence operator and it enjoys the following properties
\[
\mathcal{B}: \Big\{g\in L^{p}; \int_{\Omega}g dx=0 \Big\} \rightarrow W^{1,p}_{0}(\Omega),
\]
\[
\|\mathcal{B}(g)\|_{W^{1,p}(\Omega)} \leq C \|g\|_{L^{p}(\Omega)}.
\]
If $g$ can be written as $g=\Div h$ for a certain $h \in L^{r}$ with $h\cdot \hat{n}=0$ on $\partial\Omega$, then
\[
\|\mathcal{B}(g)\|_{L^{r}(\Omega)} \leq C \|h\|_{L^{r}(\Omega)}.
\]
We will use this operator to obtain higher integrability of $\vr_{n}$. By extending (\ref{eq:2.10}) to zero outside $\Omega$ and regularizing it, we have,
\begin{equation} 
\label{eq:4.1}
\partial_{t}b(\vr_{n})_{\epsilon}+\Div (b(\vr_{n})_{\epsilon}u)+\Big(\big[ b'(\vr_{n})\vr_{n}-b(\vr_{n})\big] \Div \vu_{n} \Big)_{\epsilon}=r_{\epsilon},
\end{equation}
where as proved in Lions\cite{Lions1998}, $r_{\epsilon} \rightarrow 0$ in $L^{2}((0,T)\times\R^{3})$. We are now ready to prove the following result.

We take a test function of the form
$$\phi_{i}=\chi(t)\mathcal{B}_{i}\Big[b(\vr_{n})_{\epsilon}- \oint_{\Omega}b(\vr_{n})_{\epsilon}dy\Big],$$ 
where
\[ \oint_{\Omega}b(\vr_{n})_{\epsilon}dy=\frac{1}{|\Omega|}\int_{\Omega}b(\vr_{n})_{\epsilon}dy, \quad \chi \in \mathcal{D}(0,T)
\]
and test it against (\ref{eq:1.10 b}). Then, with the aid of (\ref{eq:4.1}),
\begin{equation*}
 \begin{split}
 &\int^{T}_{0}\!\!\int_{\Omega} \chi \vr_{n}^{\gamma_{n}}b(\vr_{n})_{\epsilon}dxdt \\
 &= \int^{T}_{0}\!\!\int_{\Omega} \chi \vr_{n}^{\gamma_{n}} \Big[\oint_{\Omega}b(\vr_{n})_{\epsilon}dy \Big]dxdt -\int^{T}_{0}\!\!\int_{\Omega} \chi_{t}\vr_{n} \vu_{n} \cdot \mathcal{B}\Big[ b(\vr_{n})_{\epsilon}-\oint_{\Omega}b(\vr_{n})_{\epsilon}dy\Big] dxdt\\
 &+\int^{T}_{0}\!\!\int_{\Omega} \chi \vr_{n} \vu_{n}  \cdot \mathcal{B}\Big[ \big( (b^{'}(\vr_{n})\vr_{n}-b(\vr_{n}))\Div \vu_{n} \big)_{\epsilon} -\oint_{\Omega} \big( (b^{'}(\vr_{n})\vr_{n}-b(\vr_{n}))\Div \vu_{n} \big)_{\epsilon}dy\Big]dxdt\\
 & -\int^{T}_{0}\!\!\int_{\Omega} \chi \vr_{n} \vu_{n}  \cdot \mathcal{B}\Big[ r_{\epsilon} -\oint_{\Omega}r_{\epsilon}dy\Big]dxdt + \int^{T}_{0}\!\!\int_{\Omega}\chi \vr_{n}\vu_{n}  \cdot \mathcal{B}\Big[\nabla \cdot \big(b(\vr_{n})_{\epsilon} \vu_{n}  \big)\Big]dxdt\\
 &- \int^{T}_{0}\!\!\int_{\Omega}\chi \vr_{n} \vu_{ni}\vu_{nj} \partial_{i}\mathcal{B}_{j}\Big[ b(\vr_{n})_{\epsilon} -\oint_{\Omega}b(\vr_{n})_{\epsilon}dy\Big]dxdt \\
 & + \int^{T}_{0}\!\!\int_{\Omega}\chi \partial_{i}\vu_{nj}\partial_{i} \mathcal{B}_{j} \Big[ b(\vr_{n})_{\epsilon} -\oint_{\Omega}b(\vr_{n})_{\epsilon}dy\Big]dxdt\\
 & + \int^{T}_{0}\!\!\int_{\Omega}\chi \Div \vu_{n}  \Big[ b(\vr_{n})_{\epsilon} -\oint_{\Omega}b(\vr_{n})_{\epsilon}dy\Big]dxdt - \int^{T}_{0}\!\!\int_{\Omega}\chi \eta^{2}_{n}\Big[ b(\vr_{n})_{\epsilon} -\oint_{\Omega}b(\vr_{n})_{\epsilon}dy\Big]dxdt\\
 & + \int^{T}_{0}\!\!\int_{\Omega}\chi \sigma_{nij}\partial_{i} \mathcal{B}_{j} \Big[ b(\vr_{n})_{\epsilon} -\oint_{\Omega}b(\vr_{n})_{\epsilon}dy\Big]dxdt -\int^{T}_{0}\!\!\int_{\Omega}\chi \eta_{n} \Big[ b(\vr_{n})_{\epsilon} -\oint_{\Omega}b(\vr_{n})_{\epsilon}dy\Big]dxdt\\
 &=I_{1} +\cdots +I_{11}.
 \end{split}
\end{equation*}
By taking into account \eqref{e10} and the bounds of the previous sections we now estimate $I_{1}, \cdots, I_{11}$. For details, see Feireisl\cite{Feireisl2001}.\\
For $I_{1}$ we have
$$I_{1} \lesssim C(T).$$
Concerning $I_{2}$ we get
\[I_{2} \lesssim  \|\vr_{n} \vu_{n} \|_{L^{\infty}(0,T; L^{\frac{2\gamma_{n}}{\gamma_{n}+1}}(\Omega))} \|b(\vr_{n})_{\epsilon}\|_{L^{\infty}(0,T; L^{\frac{6\gamma_{n}}{5\gamma_{n}-3}}(\Omega))} \leq C(T) \|b(\vr_{n})_{\epsilon}\|_{L^{\infty}(0,T; L^{\frac{6\gamma_{n}}{5\gamma_{n}-3}}(\Omega))}.
\]
For $I_{3}$ and $I_{4}$ we get 
\begin{equation*}
 \begin{split}
 I_{3} &\lesssim  \|\vr_{n}\|_{L^{\infty}(0,T; L^{\gamma}(\Omega))} \|\nabla \vu_{n} \|^{2}_{L^{2}(\Omega\times(0,T))} \|b(\vr_{n})_{\epsilon}\|_{L^{\infty}(0,T;L^{\frac{3\gamma_{n}}{2\gamma_{n}-3}}(\Omega) )} \\
 &\leq C(T) \|b(\vr_{n})_{\epsilon}\|_{L^{\infty}(0,T;L^{\frac{3\gamma_{n}}{2\gamma_{n}-3}}(\Omega) )}.
 \end{split}
\end{equation*}
 $$I_{4} \lesssim \|\vr_{n} \vu_{n} \|_{L^{\infty}(0,T;L^{\frac{2\gamma_{n}}{\gamma_{n}+1}}(\Omega))} \|r_{\epsilon}\|_{L^{2}(\Omega\times(0,T))} \leq C(T)\|r_{\epsilon}\|_{L^{2}(\Omega\times(0,T))}.$$
We estimate now $I_{5}+I_{6}$,
\begin{equation*}
 \begin{split}
 I_{5} + I_{6} &\lesssim  \|\vr_{n}\|_{L^{\infty}(0,T; L^{\gamma_{n}}(\Omega))} \|\nabla \vu_{n} \|^{2}_{L^{2}(\Omega\times(0,T))} \|b(\vr_{n})_{\epsilon}\|_{L^{\infty}(0,T;L^{\frac{3\gamma_{n}}{2\gamma_{n}-3}}(\Omega) )} \\
 &\leq C(T) \|b(\vr_{n})_{\epsilon}\|_{L^{\infty}(0,T;L^{\frac{3\gamma_{n}}{2\gamma_{n}-3}}(\Omega) )}.
 \end{split}
\end{equation*}
For $I_{7} + I_{8}$ we have
$$I_{7} + I_{8} \lesssim \|\nabla \vu_{n} \|_{L^{2}(\Omega\times(0,T))} \|b(\vr_{n})_{\epsilon}\|_{L^{2}(\Omega\times(0,T))} \leq C(T)\|b(\vr_{n})_{\epsilon}\|_{L^{2}(\Omega\times(0,T))}.$$
and finally we get
\begin{equation*}
 \begin{split}
 & I_{9} + I_{10} + I_{11} \\
 & \lesssim \Big( \|\eta_{n}\|^{2}_{L^{2}(0,T; L^{6}(\Omega))} + \|\sigma_{n}\|_{L^{1}(0,T; L^{3}(\Omega))} + \|\eta_{n}\|_{L^{1}(0,T; L^{3}(\Omega))}\Big) \|b(\vr_{n})_{\epsilon}\|_{L^{\infty}(0,T; L^{\frac{3}{2}} (\Omega))}\\
 & \leq C(T) \|b(\vr_{n})_{\epsilon}\|_{L^{\infty}(0,T; L^{\frac{3}{2}} (\Omega))}.
 \end{split}
\end{equation*}
In sum,
\begin{equation} \nonumber
 \begin{split}
 & \int^{T}_{0}\!\!\int_{\Omega} \chi \vr_{n}^{\gamma_{n}}(b(\vr_{n}))_{\epsilon}dxdt \\
 & \leq C(T) + \|b(\vr_{n})_{\epsilon}\|_{L^{\infty}(0,T; L^{\frac{6\gamma_{n}}{5\gamma_{n}-3}}(\Omega))}+ \|b(\vr_{n})_{\epsilon}\|_{L^{\infty}(0,T; L^{\frac{3\gamma_{n}}{2\gamma_{n}-3}}(\Omega))} \\
 & + \|b(\vr_{n})_{\epsilon}\|_{L^{\infty}(0,T; L^{\frac{3}{2}}(\Omega))} + \|b(\vr_{n})_{\epsilon}\|_{L^{2}(\Omega\times(0,T))}+\|r_{\epsilon}\|_{L^{2}(\Omega\times(0,T))}.
 \end{split}
\end{equation}
By taking the limit $\epsilon \rightarrow 0$,
\begin{equation} \nonumber
 \begin{split}
 & \int^{T}_{0}\!\!\int_{\Omega} \chi \vr_{n}^{\gamma_{n}}b(\vr_{n})dxdt \\
 & \leq C(T) + \|b(\vr_{n})\|_{L^{\infty}(0,T; L^{\frac{6\gamma_{n}}{5\gamma_{n}-3}}(\Omega))}+ \|b(\vr_{n})\|_{L^{\infty}(0,T; L^{\frac{3\gamma_{n}}{2\gamma_{n}-3}}(\Omega))} \\
 & + \|b(\vr_{n})\|_{L^{\infty}(0,T; L^{\frac{3}{2}}(\Omega))} + \|b(\vr_{n})\|_{L^{2}(\Omega\times(0,T))}.
 \end{split}
\end{equation}
We approximate the function $z\mapsto z$ by a sequence of $\{b_{n}\}$ in (\ref{eq:2.10}), and approximate $\chi$ to the identity function of $(0,T)$. Then,
\begin{equation} \label{eq:4.5}
 \begin{split}
 \int^{T}_{0}\!\!\int_{\Omega} \vr_{n}^{\gamma_{n}+1}dxdt &\leq C(T) + \|\vr_{n}\|_{L^{\infty}(0,T; L^{\frac{6\gamma_{n}}{5\gamma_{n}-3}}(\Omega))}+ \|\vr_{n}\|_{L^{\infty}(0,T; L^{\frac{3\gamma_{n}}{2\gamma_{n}-3}}(\Omega))} \\
 & + \|\vr_{n}\|_{L^{\infty}(0,T; L^{\frac{3}{2}}(\Omega))} + \|\vr_{n}\|_{L^{2}(\Omega\times(0,T))}.
 \end{split}
\end{equation}
By taking into account that $\vr_{n}\in L^{\infty}(0,T;L^{1}(\Omega))$ and \eqref{e10} and sine $\gamma_{n}\to \infty$ we can always assume that $\gamma_{n}\geq N=3$ we have that the right hand side of \eqref{eq:4.5} is uniformly bounded and we can conclude that 
$$\int^{T}_{0}\!\!\int_{\Omega} \vr_{n}^{\gamma_{n}+1}dxdt \leq C(T)$$ 
which completes the proof of \eqref{e8}.\\[0.5cm]
{\bf Step 3: Convergence of  the approximating sequence $\mathbf{\{\vr_{n}, \vu_{n}, \eta_{n}, f_{n}\}}$.}\\
By using the compactness properties of  the approximating sequence $\{\vr_{n}, \vu_{n}, \eta_{n}, f_{n}\}$ stated in the Proposition \ref{prop:2.3} and the bounds of the Step 1 and Step 2 we get
\begin{equation}\nonumber
 \begin{split}
 & \vr_{n}\vu_{n} \rightarrow \vr \vu \hspace{0.2cm} \text{in} \hspace{0.2cm} L^{p}(0,T; L^{r}(\Omega)) \hspace{0.2cm} \text{for all} \hspace{0.2cm} 1\leq p<\infty, \quad 1\leq r<2, \\
 & \vu_{n} \rightarrow \vu \hspace{0.2cm} \text{in} \hspace{0.2cm} L^{p}(\Omega\times (0,T))\cap\{\vr_{n}>0\} \hspace{0.2cm} \text{for all} \hspace{0.2cm} 1\leq p<2, \\
 & \vu_{n} \rightarrow \vu \hspace{0.2cm} \text{in} \hspace{0.2cm} L^{2}(\Omega\times (0,T))\cap\{\vr_{n} \ge \delta\} \hspace{0.2cm} \text{for all} \hspace{0.2cm} \delta>0, \\
 & \vr_{n}\vu_{ni}\vu_{nj} \rightarrow \vr \vu_{i}\vu_{j} \hspace{0.2cm} \text{in} \hspace{0.2cm} L^{p}(0,T; L^{1}(\Omega)) \hspace{0.2cm} \text{for all} \hspace{0.2cm} 1\leq p<\infty,\\
&  (\vr_{n})^{\gamma_{n}}\rightharpoonup \pi, \hspace{0.2cm} \text{where $\pi\in \mathcal{M}((0,T)\times\Omega),$}\\
& \eta_{n} \rightarrow \eta,  \hspace{0.2cm} \text{in} \hspace{0.2cm} L^{2}(\Omega\times (0,T)),\\
& \sigma_{n} \rightharpoonup \sigma,  \hspace{0.2cm} \text{in} \hspace{0.2cm} L^{2}(\Omega\times (0,T)),\\
 & f_{n} \rightharpoonup f \hspace{0.2cm} \text{in} \hspace{0.2cm} L^{2}(0,T; L^{\frac{6}{5}}(\Omega\times S^{2})).
 \end{split}
\end{equation}
With the above convergence result we can pass into the weak limit in the system \eqref{4.1}-\eqref{4.4}, and we get that $\vr, \vu, \eta, f$ is a weak solution of the problem \eqref{2.14}-\eqref{2.20} provided we prove the conditions \eqref{2.18}-\eqref{2.20}. This is equivalent to the proof of 
\begin{equation}
\label{e12}
\vr\pi=\pi.
\end{equation}
Setting $s_{n}=\vr_{n}\log\vr_{n}$ and $\bar{s}=\overline{\vr\log\vr}$ its weak limit.
and using (\ref{4.1}) we get
\begin{equation} \nonumber
(\vr_n \log \vr_n)_{t} +\nabla \cdot(\vr_n \log \rho_n \vu_n)+ (\nabla \cdot \vu_n)\vr_n=0.
\end{equation}
Next, we apply the differential operator $(-\Delta)^{-1}\nabla \cdot $ to (\ref{4.2}). Then
\[
\frac{d}{dt}\Big[(-\Delta)^{-1}\nabla \cdot(\vr_n \vu_n)\Big]+(-\Delta)^{-1}\partial_{i}\partial_{j}(\vr_n {\vu_n}_{i}{\vu_n}_{j})+2\nabla \cdot \vu_n -\vr^{\gamma_{n}}-\eta_n^{2}=\]
\[(-\Delta)^{-1}\nabla \cdot (\nabla \cdot \sigma_n-\nabla \eta_n),
\]
from which we have
\[
 2\nabla \cdot \vu_n= -\frac{d}{dt}\Big[(-\Delta)^{-1}\nabla \cdot(\vr_n \vu_n)\Big]-(-\Delta)^{-1}\partial_{i}\partial_{j}(\vr_n {\vu_n}_{i} {\vu_n}_{j}) +\vr_n^{\gamma_{n}}+\eta_n^{2}
 \]
 \[ +(-\Delta)^{-1}\nabla \cdot (\nabla \cdot \sigma_n-\nabla \eta_n).\]
The last two relations yield
\begin{equation} \nonumber
 \begin{split}
 & 2\Big[(\vr_{n} \log \vr_{n})_{t}+\nabla \cdot(\vr_{n} \log \vr_{n} \vu_n)\Big]+(\vr_{n})^{\gamma_{n}+1} \\
 &= -\vr_{n}(\eta_{n})^{2}-\vr_{n}\Big[(-\Delta)^{-1}\nabla \cdot(\nabla \cdot \sigma_{n}-\nabla \eta_{n})\Big] +\frac{d}{dt}\Big[\vr_{n} (-\Delta)^{-1}\nabla \cdot(\vr_{n} \vu_{n})\Big] \\
 &+ \nabla \cdot\Big[\vr_{n} \vu_{n} (-\Delta)^{-1}\nabla \cdot(\vr_{n} \vu_{n})\Big] \\
 &+ \vr_{n} \Big[(-\Delta)^{-1}\partial_{i}\partial_{j}(\vr_{n} {\vu_n}_{i}{\vu_n}_{j}) -u^{n}\cdot \nabla (-\Delta)^{-1}\nabla \cdot(\vr_{n} \vu_{n})\Big].
 \end{split}
\end{equation}
By taking the limit $n\rightarrow \infty$ in  the last relation, we get
\begin{equation} \nonumber
 \begin{split}
 & 2\Big[\overline{s}_{t}+\nabla \cdot(\vu \overline{s})\Big]+\overline{(\vr_{n})^{\gamma_{n}+1}} \\
 &= -\vr\eta^{2}-\vr\Big[(-\Delta)^{-1}\nabla \cdot(\nabla \cdot \sigma-\nabla \eta)\Big] +\frac{d}{dt}\Big[\rho (-\Delta)^{-1}\nabla \cdot(\vr \vu)\Big] \\
 &+ \nabla \cdot\Big[\vr \vu (-\Delta)^{-1}\nabla \cdot(\vr \vu)\Big] \\
 &+ \vr \Big[(-\Delta)^{-1}\partial_{i}\partial_{j}(\vr \vu_{i}u_{j}) -u\cdot \nabla (-\Delta)^{-1}\nabla \cdot(\vr \vu)\Big],
 \end{split}
\end{equation}
where we use  Proposition \ref{prop:2.3} (iii) to pass to the limit of $\vr_{n}(\eta_{n})^{2}$. Next, we take the limit of (\ref{4.2}). By Proposition \ref{prop:2.3} (ii),
\[
\partial_{t}(\vr \vu) +\nabla \cdot (\vr \vu \otimes \vu)-\Delta \vu -\nabla (\nabla \cdot \vu)+\nabla \pi+\nabla \eta^{2}=\nabla \cdot \sigma-\nabla \eta.
\]
Let $s=\vr \log \vr$. By following the same calculations above, we obtain that
\begin{equation} \nonumber
 \begin{split}
 & 2\Big[{s}_{t}+\nabla \cdot(\vu {s})\Big]+ \vr\pi \\
 &=-\vr\eta^{2} -\vr\Big[(-\Delta)^{-1}\nabla \cdot(\nabla \cdot \sigma-\nabla \eta)\Big] +\frac{d}{dt}\Big[\vr (-\Delta)^{-1}\nabla \cdot(\vr \vu)\Big] \\
 &+ \nabla \cdot\Big[\vr \vu (-\Delta)^{-1}\nabla \cdot(\vr \vu)\Big] \\
 &+ \vr \Big[(-\Delta)^{-1}\partial_{i}\partial_{j}(\vr \vu_{i}\vu_{j}) -u\cdot \nabla (-\Delta)^{-1}\nabla \cdot(\vr \vu)\Big].
 \end{split}
\end{equation}
Comparing the last two relations, we have
\begin{equation}
\partial_{t}(\bar{s}-s)+\Div\left ((\bar{s}-s)\vu_{n}\right)=-\overline{\vr\Div\vu}+\vr_{n}\Div\vu_{n}
\label{e13}
\end{equation}
and
\begin{equation}
\partial_{t}(\bar{s}-s)+\Div\left ((\bar{s}-s)\vu_{n}\right)=\frac{1}{2}\left(\vr\pi-\overline{(\vr_{n})^{\gamma_{n}+1}}\right).
\label{e14}
\end{equation}
Next, using that
$$(\vr)^{\gamma_{n}}\rightarrow \mathbf{1}_{\{\vr=1\}}, \quad \text{a.e. in $L^{p}((0,T)\times \Omega)$},$$
which yields
$$(\vr)^{\gamma_{n}}(\vr_{n}-\vr)\rightharpoonup 0,$$
we obtain
\begin{equation}
\overline{(\vr_{n})^{\gamma_{n}+1}}-\vr\overline{(\vr_{n})^{\gamma_{n}}}=\overline{(\vr_{n})^{\gamma_{n}}(\vr_{n}-\vr)}=\overline{((\vr_{n})^{\gamma_{n}}-\vr^{\gamma_{n}})(\vr_{n}-\vr)}\geq 0.
\label{e15}
\end{equation}
From \eqref{e15} we deduce, 
\begin{equation}
\vr \pi=\vr\overline{(\vr_{n})^{\gamma_{n}}}\leq \overline{(\vr_{n})^{\gamma_{n}+1}}.\nonumber
\end{equation}
Integrating \eqref{e14} in space we get
\begin{equation}
\partial_{t}\int_{\Omega}(\bar{s}-s)dx\leq 0. \nonumber
\end{equation}
Then, since $(\bar{s}-s)|_{t=0}=0$ and by the convexity of $s$ we have $s\leq \bar{s}$ and $s=\bar{s}$.
Therefore, from \eqref{e14} we obtain
\begin{equation}
\vr \pi=\overline{(\vr_{n})^{\gamma_{n}+1}}
\label{e18}
\end{equation}
Since for any $\varepsilon>0$, there exists $\bar{n}$, such that for any $n\geq \bar{n}$ we have the property $x^{\gamma_{n}+1}\geq x^{\gamma_{n}}-\varepsilon$ and applying it to $x=\vr$ we have
\begin{equation}
(\vr_{n})^{\gamma_{n}+1}\geq (\vr_{n})^{\gamma_{n}}-\varepsilon.
\label{e19}
\end{equation}
Passing to the weak limit in \eqref{e19} and by using \eqref{e18}we end up with
$$
\vr\pi\geq \pi-\varepsilon,
$$
and, as $\varepsilon\to 0$ we conclude with
\begin{equation}
\vr\pi\geq \pi.
\label{e20}
\end{equation}
Now it remains to prove $\vr\pi\leq \pi$. Since $\vr\pi$ is not defined almost everywhere in order to give a meaning to the inequality we want to prove we denote by $\omega_{k}$ a smoothing sequence in the space and time variables  defined as follows
$$\omega_{k}=k^{4}\omega(k\cdot),$$
$$\omega\in C^{\infty}(\R^{4}),\  \omega\geq 0, \quad \int_{\R^{4}}\!\!\!\!\!\!\!\!\!-\omega\  dxdt=1,\quad spt(\omega)\in B_{1}(\R^{4}).$$
We denote by $\vr_{k}$ and $\pi_{k}$ a sequence of smooth functions defined as
$$\vr_{k}=\vr\ast\omega, \qquad \pi_{k}=\pi\ast\omega$$
and we have that
$$\vr_{k}\rightarrow \vr \quad \text{in $C([0,T];L^{p})\cap C([0,T];H^{-1})$}$$
$$\pi_{k}\rightarrow \pi\quad \text{in $W^{-1,2}\cap L^{1}(L^{q})$}$$
for any $p,q$ such that $1/p +1/q=1$.
Hence we can rewrite $(\vr-1)\pi$ as
\begin{equation}
(\vr-1)\pi=(\vr_{k}-1)\pi_{k}+(\vr-\vr_{k})\pi_{k}+(\vr-1)(\pi-\pi_{k})
\label{e21}
\end{equation}
Since $\vr_{k}\leq 1$ and sending $k$ to $\infty$ in \eqref{e21} we obtain
\begin{equation}
\vr\pi-\pi\leq 0
\label{e22}
\end{equation}
Considering together and \eqref{e20}  and \eqref{e22} we have \eqref{e12} and finally we conclude the proof of the Theorem \ref{MT2}.

\subsection{Proof of the theorem \ref{MT}}

We can observe that the proof of the Main Theorem is a consequence of the Theorem \ref{MT2}. 
The only think we have to check is that the condition \eqref{2.18} holds in the sense of distribution. This last issue is a consequence of the following lemma (for the proof we refer to \cite{LionsMasmoudi-1999}, Lemma 2.1).
\begin{lemma}
Let $\vu\in L^{2}(0,T;H^{1}_{loc}(\Omega))$ and $\vr\in L^{2}_{loc}((0,T)\times \Omega)$ satisfying 
 $$\partial_t \vr_n +\Div(\vr_n \vu_n)=0, \quad \text{in $(0,T)\times \Omega$},$$
 $$\vr(0)=\vr_{0},$$
 then the following two assertion are equivalent
 \begin{itemize}
 \item[(i)] $\Div \vu=0$, a.e. on $\{ \vr \geq 1\}$ and $0\leq \vr_{0}\leq 1$.
  \item[(ii)] $0\leq \vr\leq 1$.
 \end{itemize}
 \end{lemma}

We conclude this section with a final remark on how to obtain  the energy inequality \eqref{energy} that we require our global weak solutions have to satisfy.
By using the convergences proved in the Theorem \ref{MT2} we can pass into the weak limit in the energy inequality \eqref{4.19} and we obtain
\begin{equation}\label{energy-bis}
 \begin{split}
 &\int_{\Omega} \Big(\frac{\rho|u|^{2}}{2}+\eta^{2} +\psi \Big)(t)dx + 4\int^{t}_{0}\int_{\Omega}\int_{S^{2}} |\nabla_{\tau}\sqrt{f}|^{2}d\tau dxdt \\
 &+4 \int^{t}_{0} \int_{\Omega} \int_{S^{2}} |\nabla \sqrt{f}|^{2}d\tau dxdt +\int^{t}_{0} \int_{\Omega}\Big(|\nabla \vu|^{2} + |\Div \vu|^{2} +2|\nabla \eta|^{2}\Big)dxdt \\
 & \leq \int_{\Omega} \Big(\frac{\rho_{0}|u_{0}|^{2}}{2}+\eta^{2}_{0} +\psi_{0} \Big)dx+\liminf_{n\to \infty}\int_{\Omega}dx\frac{(\vr_{n0})^{\gamma_{n}}}{\gamma_{n}},\quad \text{a.e in $t$}.
 \end{split}
\end{equation}

Now, if we take, for any $n>2$, $\vr_{n0}=\vr_{0}$, $m_{n0}=m_{0}$ and we recall that $0\leq\vr_{0}\leq 1$, then
$$\liminf_{n\to \infty}\int_{\Omega}\frac{(\vr_{n0})^{\gamma_{n}}}{\gamma_{n}}dx=0$$
and we get the energy inequality \eqref{energy-bis}.

\section{Acknowlegments}
The work of D.D. was supported  by the Ministry of Education, University and Research (MIUR), Italy under the grant PRIN 2012- Project N. 2012L5WXHJ, \emph{Nonlinear Hyperbolic Partial Differential Equations, Dispersive and Transport Equations: theoretical and applicative aspects}. Ê K.T.  gratefully acknowledges the support  in part by the National Science Foundation under the grant DMS-1614964  and by the Simons Foundation under the Simons Fellows in Mathematics Award 267399. 
Part of this research was performed during the visit of K.T. at University of L'Aquila which was supported under the grant PRIN 2012- Project N. 2012L5WXHJ, \emph{Nonlinear Hyperbolic Partial Differential Equations, Dispersive and Transport Equations: theoretical and applicative aspects}.

\end{document}